\input amstex
\input Amstex-document.sty

\def\ni{\noindent}%
\def\ts{\thinspace}%
\def\RR{{\Bbb R}}%
\def\"#1{{\accent"7F #1\penalty10000\hskip 0pt plus 0pt}} % Umlaute! %

\scrollmode

\pageno 315

\def\shorttitle{Some Applications of Collapsing with Bounded Curvature}
\def\shortauthor{Anton Petrunin}

\topmatter %
\title\nofrills{\boldHuge Some Applications of Collapsing with Bounded Curvature}
\endtitle

\author \Large Anton Petrunin* \endauthor

\thanks *Department of Mathematics, PSU, University Park, PA
16802, USA. E-mail: \linebreak petrunin@math.psu.edu \endthanks

\abstract\nofrills \centerline{\boldnormal Abstract}

\vskip 4.5mm

{\ninepoint In my talk I will discuss the following results which were obtained in joint work with Wilderich
Tuschmann.

1. For any given numbers $m$, $C$ and $D$, the class of $m$-dimensional simply connected closed smooth manifolds
with finite second homotopy groups which admit a Riemannian metric with sectional curvature $\vert K \vert\le C$
and diameter $\le D$ contains only finitely many diffeomorphism types.

2. Given any $m$ and any $\delta>0$, there exists a positive constant $i_0=i_0(m,\delta)>0$ such that the
injectivity radius of any simply connected compact $m$-dimensional Riemannian manifold with finite second homotopy
group and Ricci curvature $Ric\ge\delta$, $K\le 1$, is bounded from below by $i_0(m,\delta)$.

I also intend to discuss Riemannian megafolds,  a generalized notion of Riemannian manifolds, and their use and
usefulness in the proof of these results.

\vskip 4.5mm

}

%\noindent {\bf Keywords and Phrases:} Cohomology, Symmetric group, Rotation group.}
\endabstract
\endtopmatter

\document

\baselineskip 4.5mm \parindent 8mm

This note is about a couple of applications and variations
of techniques developed in [CFG],
which we found jointly with W. Tuschmann.
Namely I will talk about injectivity radius estimates for
positive pinching, a generalized notion of manifolds,
and finiteness theorems for Riemannian manifolds with bounded curvature.
The purpose of this note is to give an informal explanation of ideas
in these proofs and for more details I refer the reader to [PT].

\specialhead \noindent \boldLARGE 1. Injectivity radius estimates
and megafolds\endspecialhead

Is it true that positive pinching of the sectional curvatures
of a simply connected
manifold implies some
lower positive bound for the injectivity radius,
which does not depend on the manifold?
For dimension $=3$ this was proved by Burago and Toponogov [BT].
More generally, they proved the following:

\noindent
{\bf Theorem A.} {\it  Given any $\delta>0$,
there exists a positive constant $i_0=i_0(\delta)>0$
such that the injectivity radius of any simply connected
compact $3$-dimensional Riemannian manifold with
Ricc$\ge\delta$, $K\le 1$, is bounded from below by $i_0$.}

Moreover they made a conjecture that this result
should be also true for
higher dimensions. Later on some new examples of manifolds with positively
pinched curvature were found by Alloff and Wallach, Eschenburg and Bazaikin
([AW], [E], [B])
which disprove this conjecture in general, but since then
closely related conjectures appeared on almost each list of
open problems in Riemannian geometry.
The theorem which we proved can be formulated as follows:

\noindent
{\bf Theorem B.} {\it  Given any $m$ and any $\delta>0$,
there exists a positive constant $i_0=i_0(m,\delta)>0$
such that the injectivity radius of any simply connected
compact $m$-dimensional Riemannian manifold with finite second homotopy
group and Ricc$\ge\delta$, $K\le 1$, is bounded from below by $i_0(m,\delta)$.}

Theorem B generalizes the Burago-Toponogov Theorem A
to arbitrary dimensions and is also in even dimensions interesting,
since there is no Synge theorem for positive Ricci curvature.
For sectional curvature pinching
a similar result was obtained independently by Fang and Rong [FR].

Now I will turn to one proof of this statement which is described in the
appendix of [PT] (The main part of paper contains an other proof).
This proof makes use of a generalized notion of Riemannian manifold,
which was also described by Gromov in the end of section $8_+$ of [G3],
and employs a ``tangential'' version of Gromov-Hausdorff convergence.
Here I will just give an informal analogy which describes this notion.
The formal aspects and all further details can be found in [PT].

One may think about a manifold as a set of charts and glueing mappings.
For a Riemannian manifold,
denoting the disjoint union of all charts with the pulled back metrics
by $(U,g)$, the set of all glueing maps
defines an isometric pseudo-group action by a pseudogroup $G$
on $(U,g)$.
Here is the definition of a pseudogroup action:

\noindent
{\bf Definition.}
{\it
A {\sl pseudogroup action} (or pseudogroup of transformations)
on a manifold $M$
is given by a set $G$ of pairs of the form $p=(D_p,\bar p)$, where $D_p$ is
an open subset of $M$ and $\bar p$ is a homeomorphism $D_p\to M$,
so that the following properties hold:

{(1)} $p,q\in G$ implies
$p\circ q=(\bar q^{-1}(D_p\cap \bar q(D_p)),\bar p\circ\bar q)\in G$;

{(2)} $p \in G$ implies $p^{-1}=(\bar p(D_p), \bar p^{-1})\in G$;

{(3)} $(M, id)\in G$;

{(4)} if $\bar p$ is a homeomorphism from an open set $D\i M$ into $M$ and
$D=\bigcup_\alpha D_\alpha$, where $D_\alpha$ are open sets in $M$,
then the property $(D,\bar p)\in G$ is equivalent  to
$(D_\alpha,\bar p|_{D_\alpha})\in G$ for any $\alpha$.

We call the pseudo-group action {\sl natural} if in addition
the following is true:

(i)$'$ If $(D,\bar p)\in G$ and $\bar p$ can be extended as a continuous
map to a boundary point $x\in \partial D$, then there is an element
$(D',\bar p')\in G$ such that $x\in D'$, $D\i D'$ and
$\bar p'|_{D}=\bar p$.
}

To form a manifold this action must be in addition properly discontinuous and free. If it just properly
discontinuous then we obtain an orbifold. In the case of a general (isometric!) pseudogroup action we obtain a
{\sl it Riemannian megafold} (cf. [PT]). The megafold which is obtained this way will be denoted by $({\Cal
M},g)=((U,g):G)$.

Now we come to the main notion of this section:

\noindent {\bf Definition.} {\it A sequence of Riemannian megafolds $({\Cal M}_n,g_n)$ is said to {\sl
Grothendieck-Lipschitz converge (GL-converge)} to a Riemannian megafold $({\Cal M},g)$ if there are
representations $({\Cal M}_n,g_n)=((U_n,g_n):G_n)$ and $({\Cal M},g)=((U,g):G)$ such that

(a) The $(U_{n},g_n)$ Lipschitz converge to $(U,g)$, and

(b) For some sequence $\epsilon_n\to 0$
there is a sequence of $e^{\pm\epsilon_n}$-bi-Lipschitz
homeomorphisms $h_n: (U_{n},g_n)\to (U,g)$ such that
the pseudogroup actions on $\{(U_{n},g_n)\}$
converge (with respect to the homeomorphisms $h_n$)
to a pseudogroup action on $\{(U,g)\}$.

I.e.,
for any converging sequence of elements
$p_{n_k}\in G_{n_k}(U_{n_k},{\Cal M}_{n_k})$ there exists
a sequence $p_n\in G_n$ which converges to the same local isometry on
$U$, and
the pseudogroup of all such limits, acting on $U$, coincides with
the pseudogroup action $G(U,{\Cal M})$.
}

Here are two simple examples of GL-convergence:

Consider the sequence of Riemannian manifolds $S^1_\epsilon\times \RR$,
which for $\epsilon \to 0$ Gromov-Hausdorff converge to $\RR$.
Then this sequence converges in the GL-topology to a Riemannian megafold
${\Cal M}$, which can described as follows:
It is covered by one single chart $U=\RR^2$, and the pseudogroup $G$
simply consists of all vertical shifts of $\RR^2$.
I.e., ${\Cal M}$ is nothing but
$(\RR^2: \RR)$ where $\RR$ acts by parallel translations.
(Note that $(\RR^2:\RR)\not=\RR^2/\RR$,
these megafolds even have different dimensions!)

The Berger spheres,
as they Gromov-Hausdorff collapse to $S^2$,
converge in Grothendieck-Lipschitz topology to the Riemannian megafold
$(S^2\times \RR:\RR)$.
Here $\RR$ acts by parallel shifts of $S^2\times \RR$.

Notice that a Riemannian metric on a megafold $((U,g):G)$
defines a pseudometric on the set
of $G$ orbits. In particular one has that
the diameter of a Riemannian megafold is
well defined. Now here is the basic result,
whose proof is obvious from the definitions:

\

\noindent {\bf Theorem C.} {\it The set of Riemannian $m$-manifolds  (megafolds) with bounded sectional curvature
$|K|\le 1$ and diameter $\le D$ is precompact (compact) in the \linebreak Grothendieck-Lipschitz topology. }

Now let us state some natural questions which arise from this theorem:

1. Which Riemannian megafolds can be approximated by manifolds
with bounded curvature and diameter?

Note that
the infinitesimal motions of the pseudogroup $G$ give rise to a Lie
algebra of Killing fields on a megafold $(U,g)$
from which one can recover an isometric
local action of a connected
Lie group on $(U,g)$. Let us call this group $G_o$.
Then $G_o$ is obviously an invariant of the megafold, i.e.,
does not depend on a particular representation $(U:G)$.
It follows now from [CFG] that if $({\Cal M},g)$
is a limit of Riemannian manifolds
with bounded curvature, then $G_o({\Cal M})$ must be  nilpotent.
A direct construction moreover shows that this condition is also sufficent.

(Note that since a pure $N$-structure on a simply connected manifold is given by a torus action, one also has the
following: If a megafold can be approximated by simply connected manifolds with bounded curvature, then $G_o({\Cal
M})=\RR^k$.)

2. How can one recover the Gromov-Hausdorff limit space
from a Grothendieck-Lipschitz limit?

Let ${\Cal M}=((U,g):G)$ be a GL-limit of Riemannian manifolds.
The GH-limit is the space of $G$ orbits with the induced metric,
in other words: The Gromov-Hausdorff limit is nothing but $(U,g)/G$.

Riemannian megafolds are actually not that general
objects as they might seem at first sight.
Indeed, given a Riemannian megafold $({\Cal M},g)$ we can consider
its orthonormal frame bundle $(F{\Cal M},\tilde g)$,
equipped with the induced metric.
Now consider some representation of it, say,
$(F{\Cal M},\tilde g)=((U,\tilde g):G)$.
Then the $G$ pseudogroup action is free on $U$, so that
its closure $\bar G$
also acts freely. Therefore the corresponding factor,
equipped with the induced metric,
is a Riemannian manifold $Y=(U/\bar G, \bar g)$, and
there is a Riemannian submersion
$(F{\Cal M},\tilde g)\to (U/\bar G, \bar g)$
whose fibre is $G_o/\Gamma_o$, where
$\Gamma_o$ is a dense subgroup of $G_o$
(Roughly speaking, $\Gamma_o$ is
generated by the intersections of $G_o$ and $G$).
If we assume that ${\Cal M}$ is simply connected,
then $G_o=\RR^k$ and $\Gamma_o$ is the homotopy sequence image of $\pi_2(Y)$.
In particular, the dimension of the free part of $\pi_2(Y)$ is at least $k+1$.

Notice that for Riemannian megafolds one can define the
de Rham complex just as well as for manifolds.
(In fact I am not aware of a single notion or theorem in Riemannian
geometry which does not admit a straightforward generalization
to Riemannian megafolds!)
From  the above characterization of Riemannian megafolds
it is not hard to obtain the following:

\noindent {\bf  Theorem D.} {\it Let $M_n$ be a sequence of
compact simply connected Riemannian $m$-manifolds with bounded
curvatures and diameters and $H^2_{dR}(M_n)=0$ which
Grothendieck-Lipschitz converges to a Riemannian megafold $({\Cal
M},g)$.

Then ${\Cal M}$ is either a Riemannian manifold and the manifolds
$M_n$ converge to ${\Cal M}$ in the Lipschitz sense, or
$H^1_{dR}({\Cal M})\not=0$.
}

\

It is in particular straightforward to show
that if  Ricc$({\Cal M})>0$, then $H^1_{dR}({\Cal M})=0$.
Moreover, a Grothendieck-Lipschitz limit of manifolds with uniformly
bounded sectional curvatures and Ricc$\ge \delta>0$ is a
Riemannian megafold with Ricc$\ge \delta>0$.

Now we can prove Theorem B:
Assume it is wrong.
Then we can find a collapsing sequence
of simply connected manifolds with finite $\pi_2$ and
positive Ricci-pinching, and
we obtain a megafold with $H^1_{dR}\not=0$ as a GL-limit.
Applying the Bochner formula for
$1$-forms on this megafold, we obtain a contradiction.

\specialhead \noindent \boldLARGE 2. Finiteness theorems
\endspecialhead

The following result appeared as a co-product of the theorem above,
and it came as a nice surprise.
Let me first formulate this finiteness results from [PT]:

\ni {\bf Theorem E (The $\pi_2$-Finiteness Theorem).\ts} {\it For
given $m$, $C$ and $D$, there is only a finite number of
diffeomorphism types of simply connected closed
$m$-dimen\-sio\-nal manifolds $M$ with finite second homotopy
groups which admit Riemannian metrics with sectional curvature
$|K(M)|\le C$ and diameter diam$(M)\le D$. }

\ni {\bf Theorem F (A ``classification'' of simply connected
closed manifolds).} {\it For given $m$, $C$ and $D$, there exists
a finite number of closed smooth simply connected manifolds $E_i$
with finite second homotopy groups such that any simply connected
closed $m$-dimen\-sio\-nal manifold $M$ admitting a Riemannian
metric with sectional curvature $|K(M)|\le C$ and diameter
diam$(M)\le D$ is diffeomorphic to a factor space $M=E_i/T^{k_i}$,
where $0\le k_i=b_2(M)=$ dim $E_i-m$ and $T^{k_i}$ acts freely on
$E_i$. }

Here is a short account of other finiteness results
which only require volume, curvature, and diameter bounds:
For manifolds $M$ of a given fixed dimension $m$, the conditions

\item{$\bullet$} $vol(M)\ge v>0$, $|K(M)|\le C$ and diam$(M)\le D$
imply finiteness of diffeomorphism types (Cheeger ([C]) 1970);
this conclusion continues to hold for
$vol(M)\ge v>0$, $\int_M\vert R\vert \sp {m/2}\leq C$, $|\text{Ric}_M|\le C'$,
diam$(M)\le D$ (Anderson and Cheeger ([AC1]) 1991);

\item{$\bullet$} $vol(M)\ge v>0$, $K(M)\ge C$ diam$(M)\le D$
imply finiteness of homeomorphism types
(Grove-Petersen ([GP]) 1988, Grove-Petersen-Wu ([GPW]) 1990); Perelman, ([Pe]) 1992)
(if in addition $m>4$, these conditions  imply
finiteness of diffeomorphism types)
and
Lipschitz homeomorphism types (Perelman, unpublished);

\item{$\bullet$} $K(M)\ge C$ and diam$(M)\le D$
imply a uniform bound for the total Betti number (Gromov [G1] 1981).

The $\pi_2$-Finiteness Theorem  requires
two-sided bounds on curvature, but no lower uniform volume bound.
Thus, in spirit it is somewhere between Cheeger's Finiteness
and  Gromov's Betti number Theorem.

Each of the above results has (at least)
two quite different proofs, the original one and one which
uses Alexandrov techniques.
(For Gromov's Betti number theorem we made such a proof recently,
jointly with V. Kapovich and it
turned out that one can even give an upper estimate for
the total number of critical points of a Morse function
on such a manifold, which due to the Morse inequality is a stronger condition.)
Let me now explain roughly this second way of proving of such theorems:

 I will take Cheeger's theorem as an example:
Assume it is wrong. Then there is an
infinite number of non-diffeomorphic manifolds with bounded curvature,
diameter and a  lower bound on the volume.
Then due to Gromov's compactness theorem a subsequence
of them has  a limit.
Then, due to the volume bound, this limit space has the same dimension, and
is in fact just little worse than Riemannian; it is a manifold with a smooth
structure and curvature bounded in the sense of Alexandrov.
Then one only has to prove the
stability result, i.e. one has to prove that starting
from some big number all manifolds are diffeomorphic
to the limit space. In the case of two-sided curvature
bound it is really simple, and
for just lower curvature bound it is already a hard theorem of Perelman,
but still it works along this lines.

Now for both of these proofs it is very important
to have a uniform lower positive volume bound to prevent collapsing.
In fact, if one removes this bound then it is not
hard to construct infinite sequence of non-diffeomorphic manifolds.
This holds for two-sided bounded as well as for lower curvature bound.
And if we would try to prove it the same way as before
we would get a limit space of possibly smaller dimension.
Therefore the stability result can not hold this way.

This partly explains why Theorem E looks a bit surprising, we add one topological
condition and get real finiteness result.
The proof can go along the same lines.
Take a sequence of nondiffeomorphic Riemannian manifold $(M_n,g_n)$,
by Gromov's compactness theorem we have a limit space (for some subsequence) $X$.
The sequence must collapse, otherwise the same arguments as before would work.
Since the $M_n$ are simply connected,
from [CFG] we have that collapsing takes place along some $T^k$-orbits of some $T^k$-action.

Now assume for simplicity that $X$ is a manifold and $\pi_2(M_n)=0$.
Then all $M_n$ are diffeomorphic to $T^k$ bundles over $X$.
Since the $M_n$ are simply connected so is $X$.
Therefore the diffeomorphism type of $M_n$ depends only on the Euler class
$e_n$ which in this case can be interpreted as the following mapping:
$$0=\pi_2(M_n)\to\pi_2(X)\buildrel {e_n} \over \rightarrow \pi_1(T^k)\to\pi_1(M_n)=0.$$
Therefore $e_n$ isan isomorphism between two groups
and up to automorphisms of $T^k$
all possible Euler classes $e_n$ are the same.
In particular,  for large $n$ all $M_n$
are diffeomorphic.

That is not quite a proof
since we had made quite strong assumptions on the way.
But it turns out that
the general case can be ruled out using a few already standard tricks
from [CFG] and [GK],
namely, by passing to the
frame bundles $FM_n$
and by conjugating group actions.

\specialhead \noindent \boldLARGE References \endspecialhead

\widestnumber\key{GPW1}

\ref \key AW \by S. Aloff; N. R. Wallach \pages 93--97 \paper An
infinite family of $7$-manifolds admitting positively curved
       Riemannian structures
\jour Bull. Amer. Math. Soc.
\vol 81
\yr 1975
\endref

\ref
\key B
\by Ya. V. Basaikin
\pages
\paper On a certain family of closed $13$-dimensional manifolds of
       positive curvature
\jour Siberian Mathematical Journal
\vol 37:6
\yr 1996
\endref

\ref \key BT \by Y. Burago; V. A. Toponogov \pages 881--887 \paper
On three-dimensional Riemannian spaces with curvature bounded
above \jour Matematicheskie Zametki \yr 1973 \vol 13
\endref

\ref \key C \by  J. Cheeger \pages 61--74 \paper Finiteness
theorems for Riemannian manifolds \jour  Amer. J. Math. \yr 1970
\vol 92
\endref

\ref \key CFG \by J. Cheeger; K. Fukaya; M. Gromov \pages 327--372
\paper Nilpotent structures and invariant metrics on collapsed
manifolds \jour J. A.M.S \yr 1992 \vol 5
\endref

\ref \key Es \by J.-H. Eschenburg \pages 469--480 \paper New
examples of manifolds with strictly positive curvature \jour
Invent. math. \yr 1982 \vol 66
\endref

\ref FR
\key FR
\by F. Fang; X. Rong
\paper Positive Pinching, volume and second Betti number
\jour GAFA (Geometric and functional analysis
\vol 9
\yr 1999
\endref

\ref \key G1 \by  M. Gromov \pages 179--195 \paper Curvature,
diameter and Betti numbers \jour Comment. Math. Helv. \yr 1981
\vol 56
\endref

\ref \key G2 \by  M. Gromov \pages 55--97 \paper Stability and
Pinching \jour Seminare di Geometria,
  Giornate di Topologia e geometria delle varieta.
  Universita degli Studi di Bologna
\yr 1992
\vol
\endref

\ref \key G3 \by M. Gromov \paper ``Metric Structures for Riemannian and Non-Riemannian \linebreak spaces'' \jour
Birkh\"auser, Basel, \yr 1999
\endref

%\ref
%\key GLP???????
%\by M. Gromov; J. Lafontaine; P. Pansu
%\paper Structures m\'etriques pour les vari\'et\'es riemannienes
%\jour Cedic Fernand, Paris
%\yr 1981
%\endref

\ref \key GK \by K. Grove; H. Karcher \paper How to conjugate
$C^1$-close group actions \jour Math. Z. \vol 132 \yr 1973 \pages
11--20
\endref

\ref \key GPW \by K. Grove; P. Petersen; J. Wu \paper Controlled
topology in geometry \jour Invent. Math. \vol 99 \yr 1990 \pages
205--213; Erratum: Invent. Math. 104 (1991), 221--222
\endref

\ref \key KS1 \by W. Klingenberg; T. Sakai \paper Injectivity
radius estimates for 1/4-pinched manifolds \jour Arch. Math. \vol
34 \yr 1980 \pages 371--376
\endref

\ref
\key PRT
\by  A. Petrunin; X. Rong; W. Tuschmann
\pages 699--735
\paper  Collapsing vs. Positive Pinching
\jour GAFA (Geometric and functional analysis)
\yr 1999
\vol 9
\endref

\ref
\key PT
\by  A. Petrunin and W. Tuschmann
\pages 736--774
\paper Diffeomorphism finiteness, positive pinching, and second homotopy
\jour GAFA (Geometric and functional analysis)
\yr 1999
\vol 9
\endref

\ref
\key PT2
\by  A. Petrunin and W. Tuschmann
\pages 775--788
\paper Asymptotical flatness and cone
structure at infinity
\jour Math. Ann.
\yr 2001
\vol 321
\endref

\enddocument